**Innovation Diffusion in EV Charging Location Decisions: Integrating Demand & Supply through Market Dynamics**

Xiangyong Luo
School of Sustainable Engineering and the Built Environment, Arizona State University
Tempe, AZ 85281, USA
Email: xluo70@asu.edu

Michael J. Kuby
School of Geographical Sciences and Urban Planning, Arizona State University
Tempe, AZ 85281, USA
Email: mikekuby@asu.edu

Yudai Honma
Institute of Industrial Science, The University of Tokyo
Meguro-ku, Tokyo, 153-8505, Japan
Email: yudai@iis.u-tokyo.ac.jp

Mouna Kchaou-Boujelben
College of Business and Economics, UAE University
Al Ain, United Arab Emirates
Email: mouna.boujelben@uaeu.ac.ae

Xuesong (Simon) Zhou *
School of Sustainable Engineering and the Built Environment, Arizona State University
Tempe, AZ 85281, USA
Email: xzhou74@asu.edu
(Corresponding Author)

**ABSTRACT**

This paper offers a strategic approach to Electric Vehicles (EVs) charging network planning, emphasizing the integration of demand and supply dynamics via continuous-time fluid queue models and discrete flow refueling location modeling, all in the context of innovation diffusion principles. We employ a continuous-time approximation based on Ordinary Differential Equations (ODEs) to design multi-year supply curves, a method that stands in contrast to conventional practices which often overlook inter-year transitions and ongoing processes. For medium-term charging station location planning (CSLP), we apply a flow refueling location model (FRLM) within grid-based multi-level networks, considering both multiple-path networks and capacity constraints. The grid-based network planning strategy uses a three-tier (Macro-Meso-Micro) approach for thorough EV charging station placement, with the macro-level covering entire cities, the meso-level assessing detailed EV routes and bridging the macro to micro levels, and the micro-level focusing on precise station placement for accessibility and efficiency. Our investigation into overutilization and underutilization scenarios delivers valuable insights for policymaking and cost-benefit analyses. Illustrating our approach with the example of the Chicago sketch network, we introduce an integrated demand-supply model suitable for a single region and extendable to multiple regions, thereby addressing a gap in the existing literature. Our proposed methodology focuses on EV station placement, taking into account future needs, geographical capacities, and the importance of scenario analysis, which empowers strategic resource planning for EV charging networks over extended timeframes, thus aiding the transition towards a more sustainable and efficient transportation system.

## 1. INTRODUCTION

The process of planning EV charging infrastructure is a complex task that needs a comprehensive understanding of demand and supply profiles to ensure efficient and economically viable deployment of charging stations. The model of innovation diffusion (Bass, 1969) significantly influences the demand side of this process by illustrating the adoption rates and market share accumulations over time. Bass's innovation diffusion model has been widely employed in transportation research, spanning various domains such as electric vehicles (EVs), electric vehicle charging stations (EVCS), autonomous vehicles, shared mobility services, and shared bike programs, showcasing its versatility in analyzing market adoption dynamics within the transportation sector (Massiani et al., 2015; Jensen et al., 2017; Shabanpour et al., 2018; Talebian et al., 2018; Zhang et al., 2020; Fluchs et al., 2020; Klein et al., 2020; Brdulak et al., 2021; Ghasri et al., 2021; Kavianipour et al., 2021). Presently, most infrastructure planning methods depend heavily on discrete optimization techniques that concentrate on constant and deterministic situations, using pre-determined candidate facility locations and fixed customer demand (Daskin, 1996; Shen et al., 2003; Drezner et al., 2004; Dong et al., 2014; Guo et al., 2016; Kuby et al., 2017; Kontou et al., 2019; Kavianipour et al., 2021). These models often present computational challenges, especially for large-scale applications, and usually require approximation and heuristic strategies for problem-solving in logistics.

Long-term demand and supply integration is a convoluted process involving considerable uncertainty. From the demand perspective, variables like technology innovation curves and long-term decisions play crucial roles. This intricacy is directly associated with the risk tolerance of planning agencies and decision-makers who must strike a balance between taking risks for potential gains and evading risks to protect against potential pitfalls.

In response to this complexity, our research proposes an innovative approach to model both demand and supply curves in EV charging network design. We utilize a continuous-time approximation approach within a continuous optimization context. This methodology uses the polynomial fluid queue concept (Newell, 1978; Cheng et al., 2022; Lu et al., 2022; Zhou et al., 2022), to approximate a discrete system into a continuous-time system using various mathematical models and techniques such as ODE. As demonstrated in Fig. 1, the EV charging planning process can be viewed as an ongoing progression of cumulative demand and supply, contextualized in terms of innovation diffusion. This model considers the yearly adoption rates and distinct phases such as early adopters and the early majority.

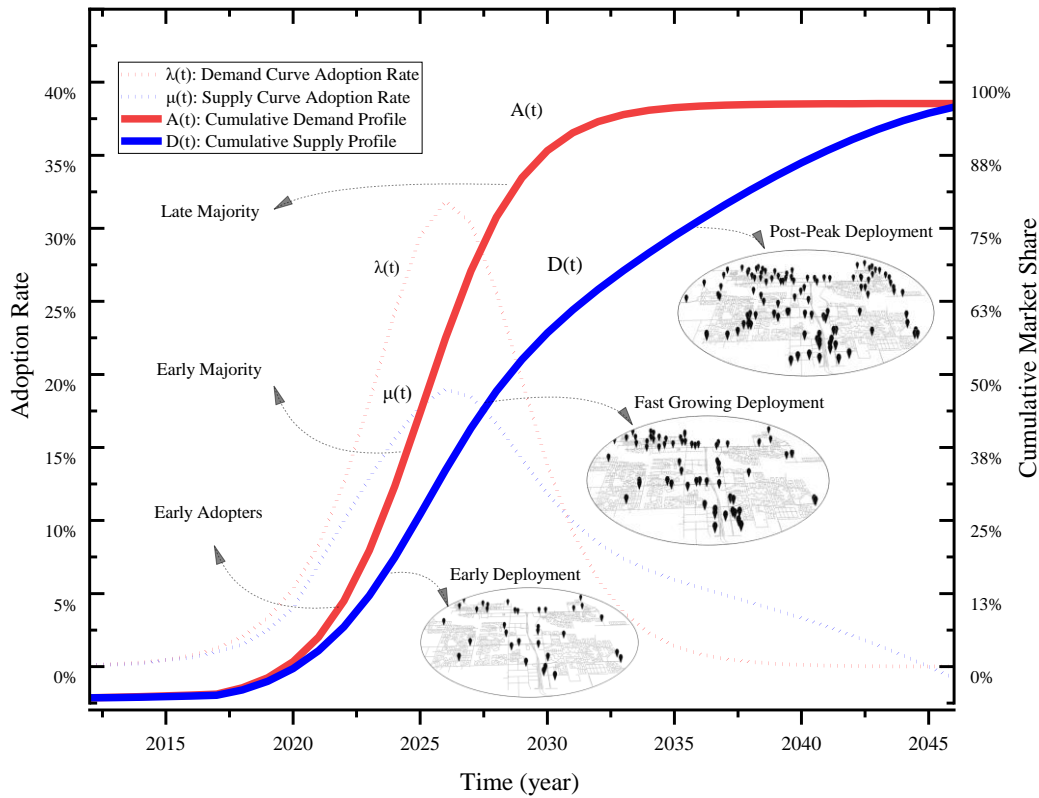

**Fig. 1** Integrative Graph Representing Long-term EV Charging Demand Market Share Prediction in Conjunction with Operational Charging Station Network Planning Strategies

The adoption of the continuous-time fluid queue approach for long-range planning is of the most significance for decision-makers in the realm of traffic flow analysis. The primary objective of this paper is to address the complexities of traffic dynamics while maintaining model simplicity and practical utility, providing crucial support for decision-makers in their EV charging network planning endeavors.

In addition, this research aims to create a simple but effective demand and supply model for EV charging planning, investigate a continuous-time strategy for long-term EV charging infrastructure, and implement this model in a single region before expanding it to multiple regions using a network framework. The study also seeks to examine infrastructure planning scenarios with both overutilization and underutilization and to develop a computational framework for predicting market adoption and saturation levels, thereby aiding resource planning decisions.

## 2. BACKGROUND

### 2.1. Innovation diffusion modeling and EV charging location demand and supply

Electric vehicle (EV) charging resource planning presents various approaches, including innovation diffusion models, discrete optimization paradigms, and continuous approximation techniques. Innovation diffusion models gained prominence in disciplines such as resource management, economics, and planning (Mansfield, 1961; Coleman et al., 1966; Bass, 1969; Rogers, 1995). They offered insights into the demand and supply dynamics of emerging technologies and were applied to various transportation planning problems, including the integration of autonomous vehicles. In contrast, discrete optimization paradigms have been employed to solve network design and facility location optimization problems (Hakimi, 1964; Church and ReVelle, 1974; Mirzain, 1985). These models were extended to account for stochastic demand features and unreliable facility services (Daskin, 1982, 1983; Dasci and Laporte, 2005; Snyder and Daskin, 2005). However, discrete optimization models often resulted in complex formulations and could be computationally demanding, especially for large-scale problems. Continuous approximation (CA) approaches and polynomial fluid queue models (Newell, 1971, 1973; Chen et al., 2020;), have been widely applied to logistics problems such as facility location, inventory management, and vehicle routing. They provided closed-form analytical structures and were often computationally efficient, particularly for large-scale practical problems. CA methods were adapted and advanced to address various emerging challenges and have been employed in the context of location, routing, and integrated supply chain and logistics studies (Ansari et al., 2018).

### 2.2. Flow Refueling Location Model

The limited driving range of alternative fuel vehicles necessitates multiple refueling stops along their routes, a consideration overlooked in the Flow Capture Location Model (FCLM) which only accounts for paths with at least one refueling facility. To address this, Kuby and Lim (2005) extended FCLM and introduced the flow refueling location model (FRLM), utilizing combinations of multiple stations to accommodate round-trip travel on longer paths. However, the immense number of facility combinations along extended routes with multiple nodes posed a significant challenge, making it impractical to generate and solve the problem even for medium-sized networks. To overcome this, Lim and Kuby (2010) devised heuristic algorithms (greedy-adding, greedy-adding with substitution, and genetic algorithms) and applied them to smaller and larger networks. There are three different formulations and solution methods that aim to avoid generating facility combinations for dealing with the range limitations. Capar and Kuby (2011) do an arc-cover path-cover (AC-PC) formulation that ensures that an entire path is covered only if each directional arc is covered. They propose a novel mixed-binary-integer programming formulation that eliminated the need for pre-generating feasible station combinations, substantially reducing the time required to solve FRLM problems while improving solution quality. A related flow-based model (Wang and Lin, 2009) considered the continuous keeping track of fuel tank or battery levels over an entire trip. MirHassani and Ebrazi (2013) made advancements in flow refueling location models (Kuby et al., 2005), introducing a flexible mixed-integer linear programming approach (expanded network approach) that significantly improved computational efficiency and outperformed the previous set-cover version in obtaining optimal solutions faster, especially for large-sized networks.

Some of other papers, in contrast, focus on model features such as covering, capacities, deviations, median approaches and so on, rather than how to model the range limitations. Mara et al., (2021) explored FRLM-related topics, including the dispersion of candidate sites on arcs, capacitated facilities, and a comparison of p-median and FRLM approaches. Wang and Wang (2010) further developed a hybrid model by combining classic set-covering model to account for nodal population demands and path demands simultaneously. New station location models extending FCLM and PMP are introduced (Hosseini et al., 2015a, 2015b, 2017a, 2017b), setting upper limits on the driving range and refueling inconvenience. Additionally, MLB-routing, a multi-path methodology based on multinomial logit, dynamically distributes packets for improved efficiency in network utilization (Honma et al., 2017). Honma and Kuby (2019) compare path-based vs. node-based models, aiming to minimize total additional travel time while feasibly covering demands with the same number of stations. Kchaou-Boujelben (2021) presented a comprehensive literature review of the charging station location problem (CSLP), critically analyzing existing studies from various perspectives, including demand representation, objective functions, side constraints, and model structures, while also considering time dependency and uncertainty in problem parameters. The investigation of electric vehicle charging station network planning and the adoption of related technologies has been the focus of several notable studies (Maia, et al, 2015; Lu et al., 2019; Gecici et al., 2022; Li et al., 2022; Song et al., 2023; Filippi et al., 2023). Incorporating multiple layers of network planning algorithms (Mahmoudi et al., 2016; Tong et al., 2017; Yang et al., 2022; Zhou et al., 2022) can significantly enhance our approach to considering both flows refueling station selection and network location deployment in large areas.

## 2.3. Limitations of existing models and research motivations

Despite numerous advancements in innovation diffusion models (Talebian, 2018; Hawkins, 2018; Zhou et al., 2018; Herrenkind, 2019; Zhang, 2020; Zhuge, 2021; Ferreira et al., 2022) and continuous approximation methods, only a few can be directly applied to real-world implementations due to specific challenges. One such challenge is accounting for overutilization and underutilization within innovation diffusion models for transportation resource planning. Additionally, the existing literature on charging station planning problems lacks comprehensive models for the continuous supply curve for charging station network planning. This emphasizes the need for a continuous-time and integrated demand and supply dynamics approach.

To bridge the methodological gaps in the existing literature, we propose a continuous-time and integrated demand and supply dynamics approach. This approach makes three significant contributions. First, it introduces a continuous-time approximation for designing supply curves for emerging technologies using Ordinary Innovation Differential Equations (OIDEs), providing efficient solution methods for large-scale transportation resource planning problems and valuable insights for system operators. Second, the proposed approach presents an analytical and parsimonious formulation with closed-form solutions, demonstrating that the original problem can be solved by addressing its corresponding revised problem, offering a more tractable solution. Additionally, the demand-supply model within a single region can be extended to multiple regions in a network design framework. This will enable multiple deployments among cities, creating a space-time network model with origins and destinations.

Within the context of this research, both long-term and medium-term CSLP perspectives hold unique meanings and applications. When discussing "long-term", we are focusing on the broad-scale, year-to-year trends in demand for EV charging stations. This is an extended view that is guided by the supply curve, providing valuable insights into the gradual shifts and trends in demand. On the other hand, the medium-term CSLP perspective pertains to the yearly or even quarterly deployment of charging stations. This level of planning is crucial for medium-term and logistical decisions, such as where and when to deploy new charging stations in response to rapidly changing demand.

In the following sections of this study, the long-term demand first provides a basis for our analysis, then we will turn our focus more explicitly to the critical aspect of medium-term charging station network planning for the deployment of charging stations.

## 3. ANALYTICAL CONTINUOUS-TIME INTEGRATED DEMAND AND SUPPLY MODEL WITH CONSIDERATION OF DIFFUSIONS OF INNOVATION

The model incorporates Newell's fluid queue model (Newell, 1978) for continuous-time demand and supply analysis, effectively representing their interplay. Additionally, the principles of Bass's innovation diffusion model (Bass, 1969) are integrated to account for the adoption of new technology over time in EV charging infrastructure planning. To facilitate our study, **Table 1** provides an overview of the relevant notations employed in this section.

**Table 1** Symbols and definitions used in extended Newell's PAQ model.

| Symbols | Definitions |
|---------|-------------|
| $t_0$ | Start time of congestion period in Newell's Model (the start time of net inconvenience of emerging technology adoption in the market) |
| $t_1$ | Time index with maximum inflow rate (maximum adoption rate of emerging technology adoption) |

| | | |
|---|---|---|
| $t_2$ | Time index with maximum queue length (time at maximum net inconvenience) | |
| $t_3$ | End time of congestion period (end time of emerging technology adoption) | |
| $\mu(t)$ | Discharge rate at time $t$ | |
| $\lambda(t)$ | Inflow rate at time $t$ | |
| $\pi(t)$ | Time-dependent net flow rate at time $t$ | |
| $D$ | Total demand during the whole peak period | |
| $P$ | Total time for congestion period | |
| $\rho$ | Net congestion curvature parameter used in polynomial form | |
| $A(t)$ | Cumulative inflow profiles at time $t$ | |
| $D(t)$ | Cumulative discharge profiles at time $t$ | |
| $Q(t)$ | Queue length at time $t$ | |
| $W(t)$ | Total delay during the congestion period at time $t$ | |

The following notations are related to Bass models

| | | |
|---|---|---|
| $p$ | Coefficient of innovation | |
| $q$ | Coefficient of imitation | |
| $m$ | Total number of innovation users in at time $t$ (EVs market size at time $t$) | |
| $t_{peak}$ | Time at the peak value of $\lambda^I(t)$ | |
| $t_{ip}$ | Inflection point time of $\lambda^I(t)$ | |
| $F(t)$ | Market share of emerging technology adoption | |
| $\lambda^I(t)$ | The demand rate of the emerging technology adoption | |
| $\pi^I(t)$ | Net inconvenience rate or generalized deviation rate of the emerging technology adoption | |
| $D^I(t)$ | Cumulative supply flow up to time t (in terms of EV charging stations) | |
| $A^I(t)$ | Cumulative demand flow up to time t (in terms of EV charging stations) | |
| $\rho^I$ | Net inconvenience curvature parameter | |
| $\mu^I(t)$ | Derived supply curve (supply adoption rate of the emerging technology adoption) | |

## 3.1 Newell's Fluid Queue Model (M1): General Demand and Supply Modeling

Newell's continuous-time fluid-based polynomial arrival queue (PAQ) model offers a valuable framework to examine the intricate relationship between demand and supply dynamics. By employing polynomial approximation for the net flow rate $\pi(t)$, we can effectively analyze continuous time arrivals and departures, along with variations in the rate of demand and supply over time. This approach allows us to consider various key elements of the model, such as the inflow rate, discharge rate, cumulative inflow, cumulative departure, queue length, and total delay time, all of which are illustrated in **Fig. 2**.

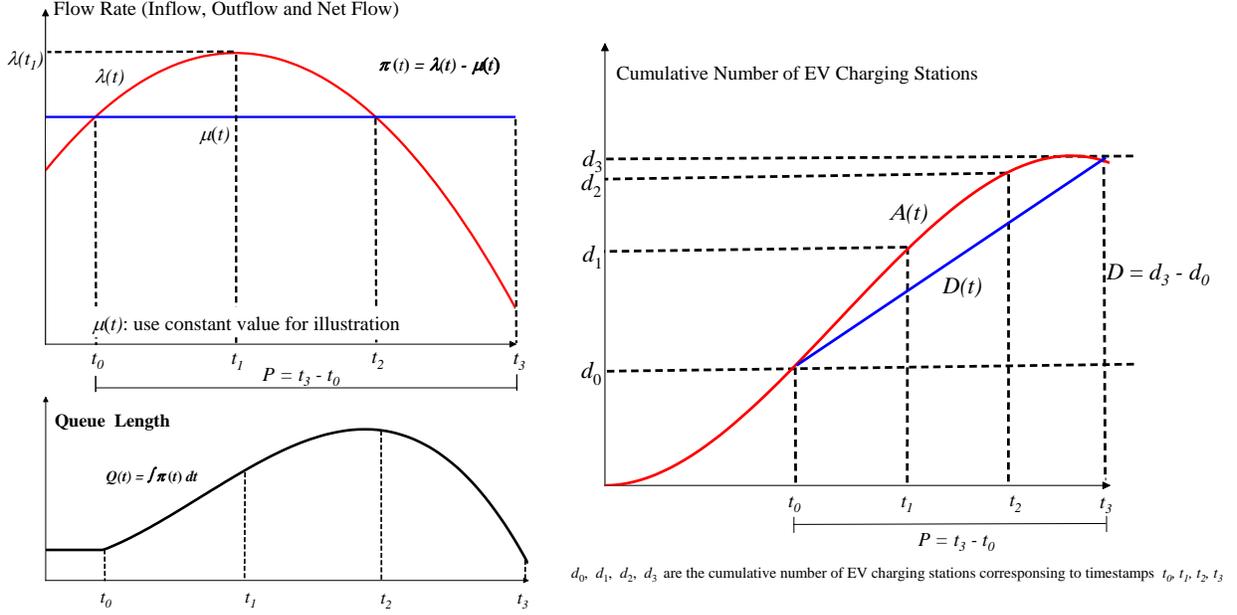

Important remark: In the foundational work on the fluid queue model, Newell used a constant $\mu(t)$ to derive an analytical form for queue length, $Q(t)$. This research examines the generalization of quadratic net flow rates $\pi(t)$ and studies the disparities in EV charging between demand and supply with an analytical form of $Q(t)$. Lu et al. (2023) recently demonstrated that discharge rates, $\mu(t)$, can be considered as time-dependent variables to calculate waiting time, $W(t)$.

**Fig. 2** Conceptual Graphical Illustration of Newell's PAQ Model for a Single Congested Period (Newell, 1982).

**Fig. 2** illustrates key parameters including the inflow rate $\lambda(t)$ at time t, start time of analysis period ($t_0$), time index with maximum inflow rate ($t_1$), time index with maximum queue length ($t_2$), end time of congestion period ($t_3$), and departure rate $\mu(t)$. The inflow rate is shown as a continuous curve, while $t_0$ and $t_3$ denote the beginning and end of the analysis period, respectively. $t_1$ represents the time at which the inflow rate reaches its maximum value during the congestion period, and $t_2$ indicates the time at which the queue length reaches its maximum value. The discharge rate $\mu(t)$ in general represents the capacity or rate of departure from the queue. The $\pi(t)$ is the net flow rate with respect to $\lambda(t) - \mu(t)$. A(t) represents the accumulated demand up to time $t$, while $D(t)$ represents the accumulated departure or dispatch from the queue up to time $t$. Queue length $Q(t)$ denotes the number of items or entities waiting in the queue at time $t$.

In the broader context of scenario planning, Newell's PAQ model can be applied to effectively analyze the diffusion and supply of emerging technologies, such as Electric Vehicles, by mapping key parameters of the model to the adoption and supply of the innovation. The inflow rate $\lambda(t)$ signifies the market demand rate of charging stations, while the departure rate $\mu(t)$ represents the market supply rate of these stations. In addition, $A(t)$ denotes the cumulative market demand at time $t$ and $D(t)$ corresponds to the cumulative market supply of charging stations. $\pi(t)$ indicates the generalized demand-supply disparity rate at time $t$ and $Q(t)$ quantifies the total generalized delay/inconvenience as a result of the dynamic interaction between demand and supply. Here are the modeling principles along with two major assumptions:

**Assumption 1**: Assume the innovation demand curve $\lambda(t)$ is given based on Bass's innovation diffusion model with parameters $p$ and $q$, and their values can be estimated from historical data or other predictions. For example, typical values for a 20-year horizon might be derived from past adoption rates of similar technologies or market trends.

**Assumption 2**: Utilize Assumption 1 to assume that the net flow rate $\pi(t)$ follows a polynomial trend (Ferreira et al., 2022). This simplification allows us to model the net flow rate over time using a quadratic function, as shown in **Eq. (2)**. Luo et al. (2023) offer a comprehensive account of the intricate data processing utilizing Uber's driver-passenger adoption case, bolstering the foundation of this assumption.

The detailed derivation of Newell's fluid queue for a closed-form solution will be given below.

Define $\pi(t)$ for a more general case with time-dependent $\lambda(t)$ and $\mu(t)$:

$$\pi(t) = \lambda(t) - \mu(t) \tag{1}$$

Based on **Eq. (1)**, time-dependent net flow can be approximated from $t_0$ to $t_3$ using a polynomial functional form. According to the **Assumption 2**, the quadratic form for net flow shows in Eq. (2).

$$\pi(t) = \rho(t - t_0)(t - t_2) \tag{2}$$

As $t_1$ represents the time with the maximum inflow rate (see Fig. 2). Newell firstly assumed that the inflow rate about time $t_1$ could be approximated by the quadratic Taylor expansion (Newell, 1986):

$$\lambda(t) = \lambda(t_1) + \lambda'(t_1) \cdot (t - t_1) + \frac{1}{2}\lambda''(t_1) \cdot (t - t_1)^2 \tag{3}$$

Since $\lambda'(t_1) = 0$, let $\rho = -\frac{1}{2}\lambda''(t_1)$, then **Eq. (3)** can be transformed to

$$\lambda(t) = \lambda(t_1) - \rho \cdot (t - t_1)^2 \tag{4}$$

Since $\rho$ describes the curvature or shape of the time-dependent inflow arrival rates. By applying the dispatch rate within a single queue duration, the queue discharge rate or service rate $\mu(t)$, where $\mu(t) = \lambda(t_0) = \lambda(t_2)$, as shown in **Fig. 2**, can be estimated in terms of **Eq. (4)**:

$$\lambda(t_1) - \rho \cdot (t_0 - t_1)^2 = \lambda(t_1) - \rho \cdot (t_2 - t_1)^2 \tag{5}$$

Now we can write $\lambda(t) - \mu(t)$ as net flow rate $\pi(t)$ by a factored form:

$$\pi(t) = \lambda(t) - \mu(t) = \rho \cdot (t - t_0) \cdot (t - t_2) \tag{6}$$

Based on the equation, the virtual queue length at time $t$ equals $A(t) - D(t)$ can be obtained:

$$Q(t) = A(t) - D(t) = \int_{t_0}^{t}[\lambda(\tau) - \mu(t)]\,d\tau \tag{7}$$

The inflow curvature parameter $\rho$ determines the shape of queueing function about $t$, further the shape of derived speed and delay profiles. The total delay between the time $t_0$ and $t_3$ can also be calculated by the area between $A(t)$ and $D(t)$ in Fig. 2 through integrating **Eq. (7)**, **Eq. (8)** is the derived total delay time:

$$W = \int_{t0}^{t}Q(t)dt = \frac{\rho}{36}(t_3 - t_0)^4 = \frac{9[\lambda(t_1) - \mu(t_1)]^2}{4\rho} \tag{8}$$

### 3.2 Innovation Diffusion Model as Demand Side Modeling (M2)

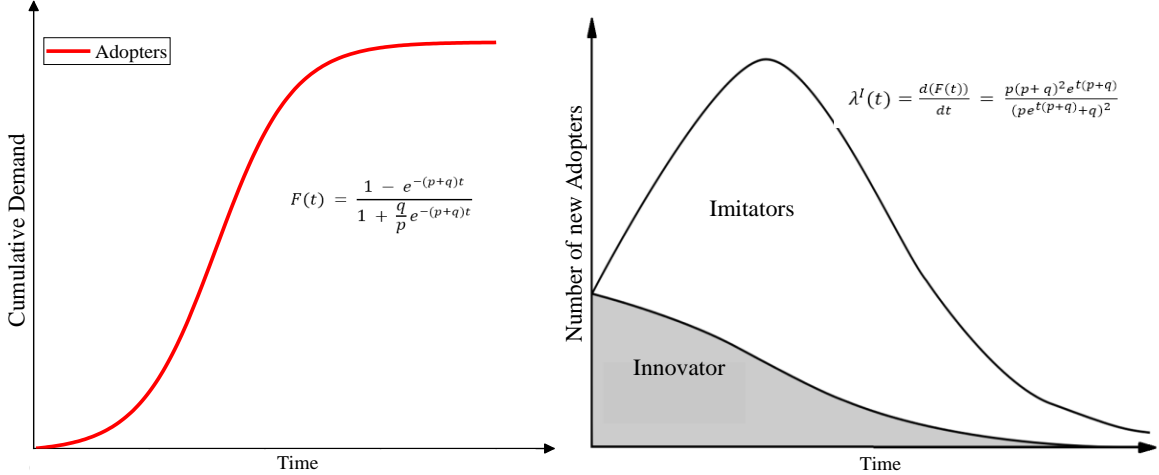

**Fig. 3** Graphical Illustration of the Bass's Innovation Diffusion Model

The innovation diffusion model will be denoting a demand curve $\lambda(t)$ in **Eq. (1)**. The innovation diffusion model delineates the market acceptance of new technologies from the demand side, taking into consideration the temporal dimension (**Fig. 3**). This model offers several advantages, including its ability to capture long-term trends and account for nonlinear relationships. In this paper, the Bass model is employed to estimate the market demand $\lambda^I(t)$ of the new technology. Notably, F(t) represents the installed base fraction, while $\lambda^I(t)$ denotes the change in the installed base fraction, defined as $\lambda^I(t) = \frac{dF(t)}{d(t)}$ (Bass, 1969). The coefficient of innovation is denoted as $p$, and the coefficient of imitation as $q$. The derivation of the innovation diffusion model is shown in **Eq. (9)**.

$$\frac{\lambda^I(t)}{1 - F(t)} = p + qF(t) \tag{9}$$

The essential ordinary differential equation is:

$$\lambda^I(t) = \frac{dF(t)}{d(t)} = p(1 - F(t)) + q(1 - F(t))F(t) \tag{10}$$

with
$$F(t) = m\frac{1 - e^{-(p+q)t}}{1 + \frac{q}{p}e^{-(p+q)t}} \tag{11}$$

which leads to

$$\lambda^I(t) = \frac{d(F(t))}{dt} = m\frac{p(p+q)^2 e^{t(p+q)}}{(pe^{t(p+q)}+q)^2} \tag{12}$$

The peak time for $\lambda^I(t)$ (Orbach, 2016; Brdulak et al., 2021) can be analytically obtained as Eq. (13).

$$t_{peak} = \frac{1}{p+q}ln(\frac{q}{p}) \tag{13}$$

The inflection point for $\lambda^I(t)$ can be derived in **Eq. (14)**.

$$t_{ip} = t_{peak} \pm \frac{ln(2+\sqrt{3})}{p+q} \tag{14}$$

The Bass's Innovation Diffusion model operates within a specific scope, concentrating primarily on the intricate mechanisms of innovation spread, notably through the lens of supply growth and temporal delays. This targeted approach allows the model to dissect the fundamental drivers of adoption while deliberately omitting complex factors such as intricate human behavioral dynamics or the influences of early adopter evaluations. By focusing on these specific aspects, the Bass's Innovation Diffusion model offers a refined framework for understanding and predicting the progression of innovations within a controlled context, enabling a more precise analysis of diffusion patterns.

### 3.3 Supply Curve Derivation (M3) Considering Underutilization and Overutilization

The supply curve μ(t) in the temporal dimension is scrutinized from the perspective of supply measurement, building upon the extension of Newell's fluid queue model. One notable demand-side model among these approaches is Bass's innovation diffusion model. It assumes that the user's total inconvenience is denoted as $\pi^I(t)$ and postulates a nonlinear temporal relationship for the total inconvenience. The subsequent steps outline the derivation process: The supply rate in the time horizon donated in **Eq. (15)**.

$$\mu^I(t) = \lambda^I(t) - \pi^I(t) \tag{15}$$

Combining **Eqs. (2) and (12)**, the derived equation in Eq. (16):

$$\mu^I(t) = \frac{p(p+q)^2 e^{t(p+q)}}{(pe^{t(p+q)}+q)^2} - \rho^I(t-t_0)(t-t_2) \tag{16}$$

To gain a clear understanding of the impact of each decision variable on the model's results (**Eq. (16)**), conducting sensitivity analysis is essential. The outcomes of this analysis have been visually depicted in the form of a detailed graph, providing insightful implications.

Fig. 4 and Fig. 5 show the results of sensitivity analysis for *p* and *q*. The innovation diffusion model encompasses two pivotal parameters that shed light on distinct facets of the process: the innovation of emerging technology (*p*) and the imitation dynamics (*q*) within the market. Regarding the innovation aspect (Fig. 4), represented by parameter *p*, as the value of *p* increases, the corresponding $t_1$ diminishes, the supply curve adoption rate and cumulative supply demonstrate a prolonged peak time. Consequently, a larger value of *p* effectively narrows the gap between demand and supply, fostering a more harmonized market equilibrium. Shifting the focus to the imitation aspect, represented by parameter *q*, Fig. 5 reveals that a higher value of *q* leads to an earlier occurrence of the peak supply curve change rate while simultaneously elevating the magnitude of the peak supply curve change rate.

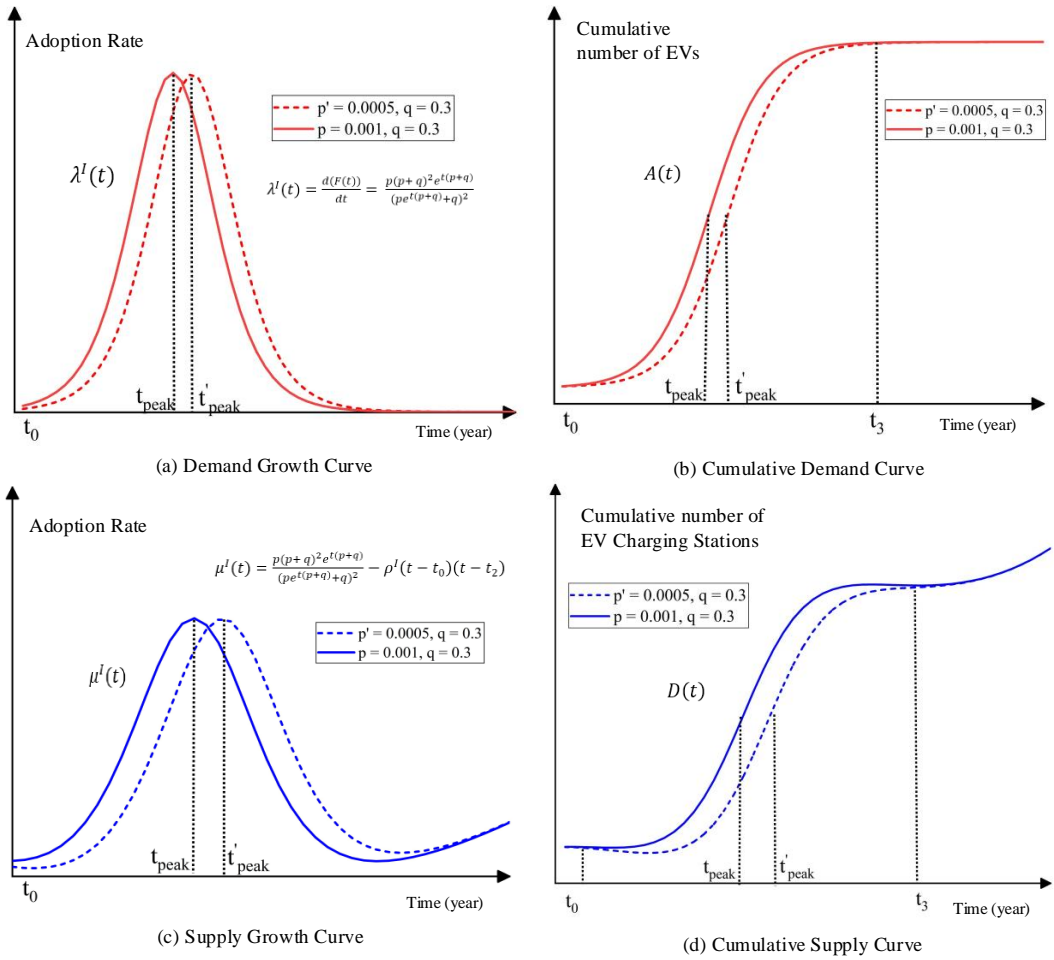

**Fig. 4** Sensitivity analysis of innovation parameter $p$ ($q$ and $\boldsymbol{\rho^I}$ are controlled)

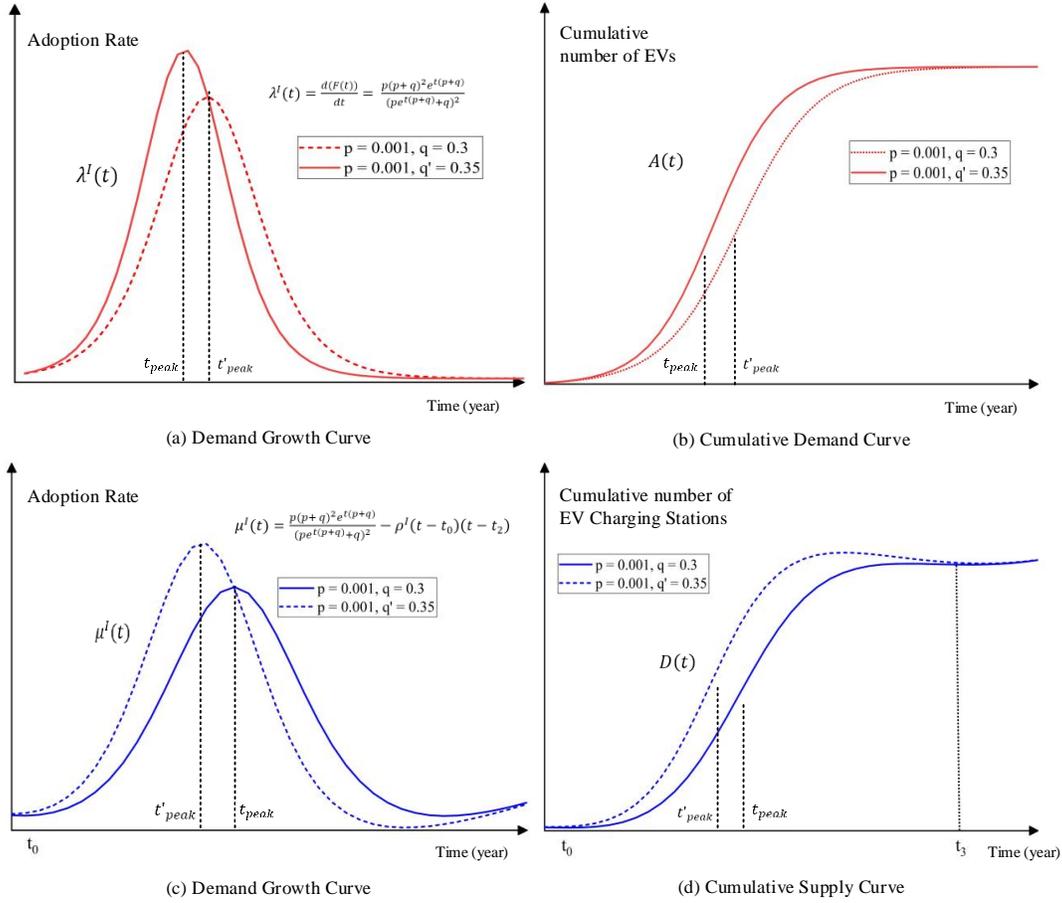

**Fig. 5** Sensitivity analysis of imitation parameter $q$ (with $p$ and $\boldsymbol{\rho^I}$ are controlled)

A thorough examination of both parameters, p, and q bestows invaluable insights into the intricate dynamics of the innovation diffusion model. Higher values of $p$ facilitate an expedited adoption of emerging technology, while concurrently harmonizing the temporal alignment of demand and supply. Larger values of $q$ contribute to earlier and intensified peak values in both demand and supply rates, thereby enriching our comprehension of the innovation diffusion process within the market landscape. These discerning findings hold paramount significance for informing strategic planning and adept management of emerging technologies, thereby ensuring their successful integration and market penetration.

Fig. 6 demonstrates a clear positive correlation between $\rho^I$ and both the final supply curve and cumulative supply curve. This correlation indicates that higher values of $\rho^I$ lead to an increased supply rate and cumulative supply. Consequently, a narrower gap between demand and supply emerges, reflecting a more balanced state. These significant findings underscore the importance of $\rho^I$ in achieving system stability and equilibrium.

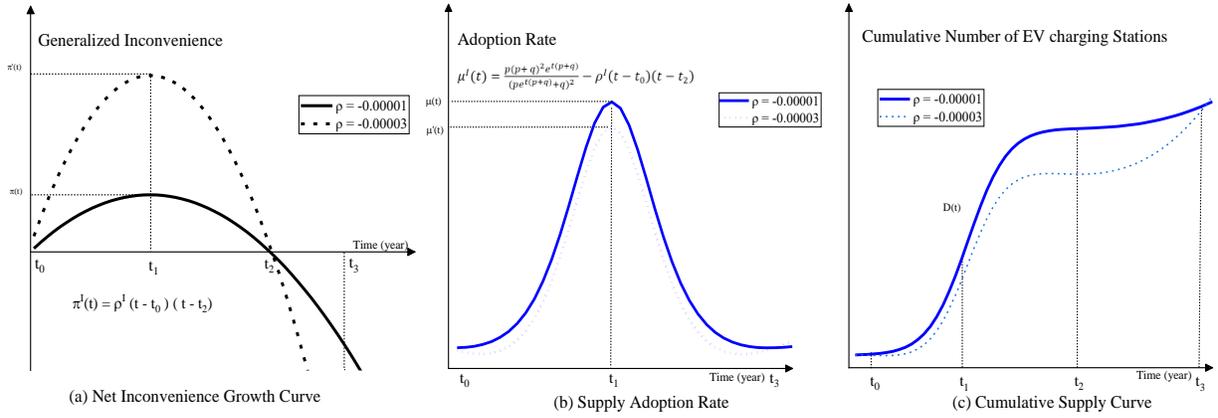

**Fig. 6** Sensitivity analysis of net inconvenience parameter $\rho^I$ (with $p$ and $q$ controlled)

According to the sensitivity analysis from Fig. 4-6, we can find that by changing these parameters, the dynamics between $\lambda^I(t)$ and $\mu^I(t)$ will represent two scenarios between the dynamic of demand and supply, which are overutilization and underutilization, as shown in Fig. 7.

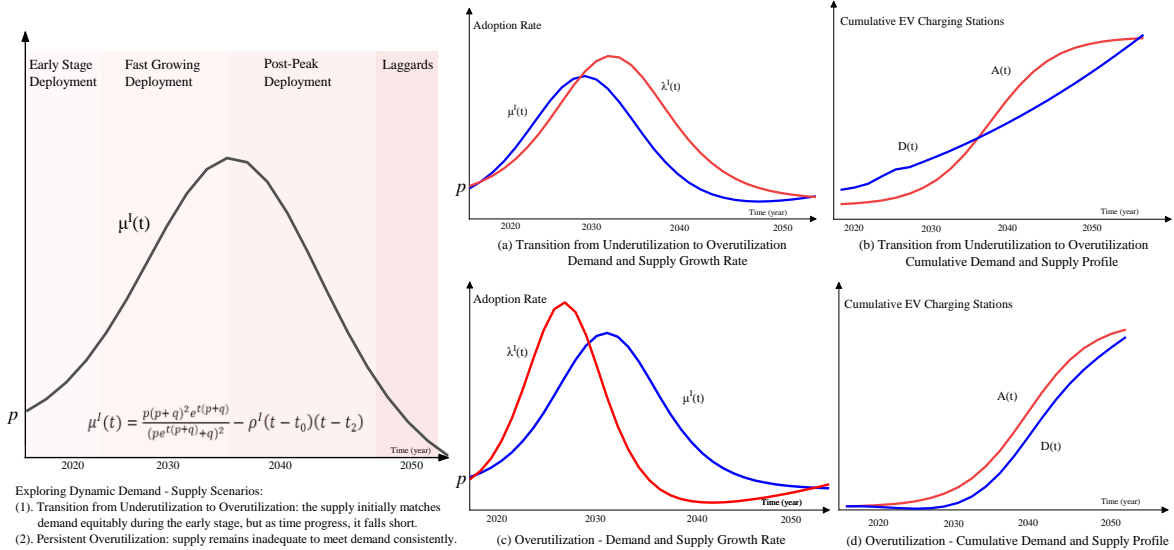

**Fig. 7** Evolution ($\mu^I(t)$) of Emerging Technology Adoption Modeled with Quadratic Form, along with Illustration of Overutilization and Underutilization Effects on Demand and Supply Dynamics

Fig. 7 portrays the derived supply curve and the general interplay between inconvenience over time and the dynamics of supply and demand. The proposed continuous approximation model extends the traditional transit dispatching problem by considering the underutilization/overutilization of emerging technology supply curves for long-term planning (Fig. 7). This integrated approach allows for more effective and sustainable transportation resource planning, ensuring that both demand and supply sides are adequately addressed. The study formulates a continuous-time approximation to design supply curves for diffusions of innovations by adapting the principle of ODEs. By considering both overutilization and underutilization scenarios, the findings will be crucial for policymakers to develop infrastructure plans and conduct cost-benefit analyses.

### 3.4 Integrating Decision Horizons in EV Planning for Managing Demand-Supply Deviations and Queueing Delays

In the context of EV planning, managing the parameter $\rho^I$ in $\pi^I(t)$ associated with **Eq. (1)** and **Eq. (16)** serves as a critical control mechanism for decision-makers. It enables them to delineate the potential deviations in the cumulative demand and supply. These deviations, in turn, may result in varying degrees of queueing at physical EV charging stations. This aspect of planning becomes more complex as it necessitates consideration of another significant coupling constraint, $D^I(t) = p^{max}$(maximum number of charging stations demanded in the region, see in Table 2), which is intrinsically tied to both long-term and CSLP decision horizons.

This intersection of medium-term supply side planning and long-term market adoption presents a unique set of challenges to planning agencies and decision-makers. Not only does it add a layer of complexity to their decision-making process, but it also exposes them to a spectrum of risks. They are compelled to strategically maneuver between risk-taking and risk aversion - a delicate balance that calls for a keen understanding of potential rewards and potential pitfalls. This risk behavior influences their approach to EV planning, highlighting the intricate connection between the two decision horizons. The decision-makers' ability to manage these risks while making efficient use of resources is an essential factor in successful EV network design and innovation diffusion.

In addition to the existing decision mechanisms, integrating a feedback process is imperative to bridge the gap between the supply side location model and the long-term continuous time model. This feedback mechanism feeds the queueing delay data from the supply side model back into the long-term model, thereby providing a dynamic framework for decision-makers to adjust their strategies based on real-time changes in queueing delay. Notably, this feedback process also helps in managing the deviation parameter, $\pi^I(t)$, effectively. The deviation information, which accounts for discrepancies between demand and supply over time, can be adjusted in line with the insights gained from the queueing delays in the supply side model.

## 4. IMMEDIATE STRATEGIES FOR EV CHARGING STATION NETWORK PLANNING (M4)

Given $p^{max} = D^I(t)$ from the long-term prediction in Section 3, this section delves into the intricacies of supply side EV charging station network planning, aiming to optimize deployment strategies and enhance accessibility for EV users. **Table 2** presents a comprehensive list of notations for key parameters that will guide our planning process and contribute to a more sustainable urban environment. The overall decision support steps are shown in Fig. 8.

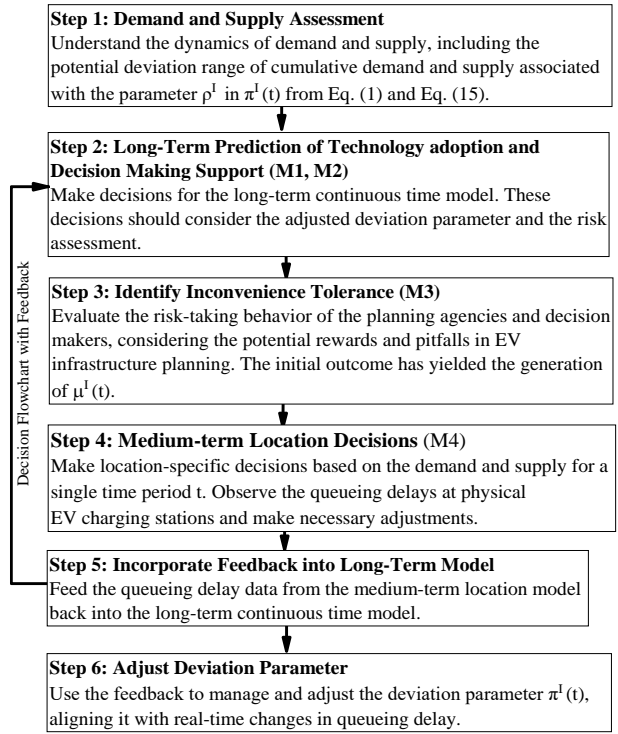

**Fig. 8** Sequential Decision Flowchart for Decision-Makers in Optimizing Electric Vehicle Charging Station Locations

The feedback mechanism outlined in Step 5 serves as a pivotal junction for decision-makers, catalyzing the selection of charging station locations. Multiple key features embedded within this step establish a robust feedback loop that significantly influences the dynamics of the overall process. Notably, Step 5 offers valuable insights, such as the Demand-to-Capacity (D/C) ratio (Zhou et al., 2022) within the network. As this ratio grows, signaling an increasing demand relative to capacity, it feeds back to Step 2, intensifying the imperative for charging station deployment.

Another crucial indicator, the Electric Vehicle / Charging Station (EV/CS) ratio (McDonald, 2019), emerges as a powerful feedback feature. This ratio, reflecting the balance between EVs and available charging stations, holds the potential to amplify inconvenience. This influence extends back to Step 2, affecting parameters $\rho^l$ and $t_2$ in the net inconvenience model. Additionally, the market budget or planning budget wields its influence across various planning stages, predominantly impacting the operational location decisions in Step 4. This financial parameter becomes instrumental in steering the progression of these stages. The intricate interplay between Step 5 and other steps, particularly the reverberations felt across Step 2 and Step 4, reflects the agility and adaptability of our integrated approach. As these feedback features dynamically influence decision-making, the model demonstrates an intricate understanding of the evolving EV landscape, infusing pragmatic realism into every stage.

### 4.1 Notation
The notations used in this section are presented in Table 2.

**Table 2** Symbols and definitions used in FRLM.

| Symbol | Definition |
|---|---|
| $i$ | Index of EV charging station in the region of interest |
| $j$ | Index of EV charging station in the region of interest |
| $k$ | Index of EV driving path |
| $Q$ | Set of all EV traveling paths in the region of interest |
| $N$ | Set of all EV traveling path-related EV charging stations |
| $N^k$ | Set of candidates charging stations for path $k$, $k \in Q$ |
| $A^k$ | Set of arcs belonging to the expanded network for path $k$, $k \in Q$ |
| $c_i$ | Generalized cost of building EV charging station at station index number $i$, $i \in N$ |
| $Cap_i$ | Maximum number of EVs could be served at a given planning horizon for charging station index number $i$ in the region, $i \in N$ |
| $Flow_k$ | Flows on traveling path $k$ at a given planning horizon, $k \in Q$ |
| $d(i, j)$ | Distance between two candidate stations $(i, j)$, $i, j \in N$ |
| $M$ | Large number in the big M method for the optimization model |
| $ord_k(i)$ | Ordering index of station $i$ in the path k, $i \in N, k \in Q$ |
| $R$ | Given EV driving range (km), the maximum driving distance with a full battery |
| $p^{max}$ | Maximum number of charging stations demanded in the region |
| **Decision Variables** | |
| $y_i$ | 1 if the EV charging station is located at node or candidate site $i$, $i \in N$; 0 otherwise |
| $x_{ij}^k$ | Flow volume on an arc $(i, j)$ of path $k$, where arc $(i, j) \in A^k$, $i, j \in N, k \in Q$ |

### 4.2 Multi-Paths EV Charging Location Selection Model

We start with the simplest path illustration (MirHassani and Ebrazi, 2013), four candidate locations within a path and the EV driving range set to 100 km, see Fig. 9. The goal is to indicate all valid combinations of stations on the path with two rules: (1) distance between two stations must be less than the driving range. (2) the amount of fuel in the tank at OD (Kuby and Lim, 2005) assumed a half-full tank in the origin. Algorithm 1 provide the data processing algorithm to construct the EV charging network for FRLM.

---

**Algorithm 1: Algorithmic Approach to Constructing an EV Charging Network for Optimization using the simplest path from Fig. 9 (a) as an example for illustration.**

---

Step 1: Initialize the node set $N^k$ and arc set $A^k$ for path k. $N^k$ : {A, B, C, D}; $A^k$ : { }. Nodes A, B, C and D are candidate EV charging stations, $k$ is the index for the selected EV charging path.

Step 2: Add two nodes o and d to the node set $N^k$, where $o$ represents the origin and $d$ refers to the destination of the traveling path. The updated node set $N^k$ : {o, A, B, C, D, d}

Step 3: Connect node o from the node set $N^k$ to any other candidate node in $N^k$ that are within EV driving range. Combined with rule (2), select all candidate nodes that are within half of the driving range to the first candidate station (A is the first candidate charging station in the selected path, see Fig. 9(a)).

For any $i \in N^k$ with $d_k(A, i) \leq \frac{R}{2}$, add arc $(o, i)$ to arc set $A^k$

Two EV driving arcs are added to the arc set and the updated arc set $A^k$: {(o, A), (o, B)}

Step 4: Connect node d from node set $N^k$ to any other node in $N^q$ that are within EV driving range. Combine with rule number (2), we select all nodes with half of the driving range to the last candidate station (D is the last candidate charging station in the selected path, see Fig. 9(a)).

For any $i \in N^k$ with $d_k(i, D) \leq \frac{R}{2}$, add arc $(i, d)$ to arc set $A^k$

Two EV driving arcs were added to set $A^k$: {(o, A), (o, B), (C, d), (D, d)}

Step 5: Connect each node $i$ of path $k$ to any other node $j$ based on two principles, (a) if the ordering index of station $i$ is less than the ordering index of station $j$. (b) the vehicle can start from node $i$ with a full tank and reach to node $j$.

for any $i, j \in N^k$ with $(ord_k(i) < ord_k(j))$ and $(d_k(i, j) \leq R)$, add arc $(i, j)$ to arc set $A^k$

Available arcs are added to the arc set $A^k$: {(o, A), (o, B), (A, B), (A, C), (B, C), (B, D), (C, D), (C, d), (D, d)}

Step 6: Finalize the node set and arc set for the optimization approach.

---

To identify the most cost-effective charging stations for EVs along multiple travel paths, we consider three station selection scenarios: (1) candidate stations have no capacity constraints, (2) candidate stations with exact capacity constraints at specific times, and (3) candidate stations with flexible capacity allowing oversaturation. Fig. 9 illustrates EV traveling flows on a single path with exact capacity constraints, the total flow on the given path (o-A-B-C-D-d) is 80 vehicles per week, and each candidate station has a capacity of 60 vehicles/week. The optimal station selections are A, B, C and D. Specifically, 60 vehicles per week will charge at stations B and D, while 20 vehicles/week will be assigned to station A and C. This approach facilitates effective and efficient station selection, considering capacity constraints for EV charging.

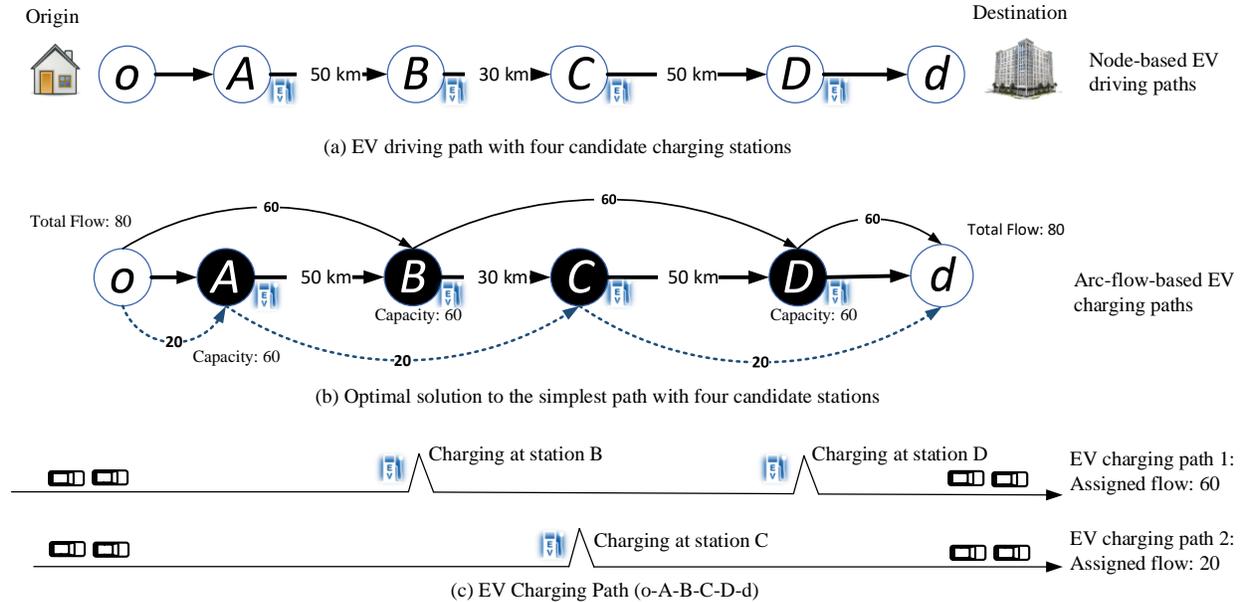

(a) EV driving path with four candidate charging stations

(b) Optimal solution to the simplest path with four candidate stations

(c) EV Charging Path (o-A-B-C-D-d)

**Fig. 9** Illustration of EV Charging Station Selection and Multiple EV Charging Paths for Single Driving Path (o-A-B-C-D-d) Considering Station Capacity Constraints

Next, we delve into the FRLM designed specifically for the selection of EV charging stations along multiple paths. The model considers various factors, including capacity constraints and flow dynamics, and EV driving ranges, to optimize the charging station deployment process efficiently. The mathematical relationship between long-term EV market planning and the location of operational facilities is expressed as follows: The total charging demand at a given time is denoted as $\sum_k Flow_k = A^I(t)$, with the total planned number of EV charging stations or points as $D^I(t)$.

**Problem statement of Model M4**. The problem at hand revolves around devising an efficient algorithm for charging stations, considering EV charging stations with and without capacity constraints, as depicted in the earlier illustrations. The primary objective is to determine the optimal charging locations for a large network with multiple EV driving paths while minimizing overall cost (Eq. (21)). This entails carefully selecting the minimum number of charging stations needed to adequately fulfill the charging requirements along each path (Eqs (17), (18)). The challenge lies in striking a balance in the flow of EVs on different paths while ensuring that the chosen charging stations effectively cater to the charging demands (Eqs. (19), (20)). By addressing this problem, we aim to establish a cost-effective and well-organized charging infrastructure that can support the widespread adoption of electric vehicles in our transportation network.

**The mass balance constraints**. For each EV charging station candidate site $i$ and each path $k$ through that node, ensure that the outflow minus inflow is equal to the total flow of the respective path (outflow - inflow). This ensures that any selected charging stations effectively handle the flow of electric vehicles along each path, maintaining a balanced charging process and meeting the overall charging demands in the network.

$$\sum_{j\,|\,(i,j)\in A^k} x_{ij}^k - \sum_{j\,|\,(i,j)\in A^k} x_{ji}^k = \begin{cases} Flow_k, & i = o \\ 0, & i \neq (o,d) \\ -Flow_k, & i = d \end{cases}, \forall k \in Q, \forall i \in N^k \quad (17)$$

**Station selection and capacity constraint**. The station selection constraint dictates that the flow can only pass through a node if a charging station is present at that specific location. This constraint ensures that charging stations are strategically positioned in the network, allowing the flow of electric vehicles only through nodes equipped with charging facilities. By adhering to $Cap_i$, the station capacity constraint becomes the smaller value. The optimization model thus endeavors to establish a harmonized and equitable charging process, efficiently regulating the flow of electric vehicles at each station within the transportation network.

$$\sum_{k \in Q} \sum_{j\,|\,(i,j)\in A^k} x_{ji}^k \leq Cap_i\, y_i, \forall i \in N^k \quad (18)$$

**Construction budget constraint**. It is imperative to ensure that the total construction cost does not surpass the predetermined financial budget denoted by $p^{max}$ and the optimization model aims to develop a feasible and cost-effective transportation network that maximizes the utilization of available resources while adhering to budgetary limitations.

$$\sum_{i \in N} y_i \leq p^{max} = D^I(t), for\ predicted\ time\ t, D^I(t)\ is\ the\ cumulative\ supply\ of\ EVCS \quad (19)$$

**Flow non-negativity constraint**. The flow non-negativity constraint is established to ensure that the flow on each arc within the network is maintained at a value greater than or equal to zero. This constraint guarantees that the flow of electric vehicles or any other relevant quantity remains non-negative, preventing negative or backward flows in the system.

$$x_{ij}^k \geq 0, for\ \forall (i,j) \in A^k, \forall k \in Q \quad (20)$$

**Objective function with hard capacity at the charging station.** The objective function minimizes the overall construction cost, given the predetermined budget input and based on capacity restrictions and supply curve predictions over time. The total number of EV charging stations is determined dynamically, considering the interplay between demand and supply as discussed in Section 3.4. By addressing this objective, the model aims to create an optimal charging infrastructure that meets the evolving supply for electric vehicles, facilitating the widespread adoption of electric vehicles in a cost-effective manner.

$$\mathbf{Z} = \sum_{i \in N} c_i y_i \quad (21)$$

We distinguish hard capacity constraints (Eq. 21), which pertain to aggregated capacity considerations (Klose et al., 2007; Wu et al., 2006). We consider the soft capacity constraints influenced by multiple factors as well. Kchaou-Boujelben (2021)'s work provides an overview of distinct capacity constraints relevant to charging station location decisions. To elaborate, we categorize the soft capacity constraints into EVs queuing theory (Dong et al., 2016), cost of average waiting time (Ghamami et al., 2016), chance constraints (Yang et al., 2018; Xie et al., 2018) and robust formulations (Upchruch et al., 2009; Mak et al., 2013; Sun et al., 2019). We integrate these concepts with Newell's PAQ model, leading us to redefine the objective function. This redefined function incorporates time/delay costs within the queuing system, effectively combining queuing theory with congestion management. The updated objective function, as expressed in **Eq. (22)**, considers station-specific arrival and departure patterns, converting delays into associated costs. $w_i$ quantifies the delay at station $i$ and VOT is the value of time parameter.

$$Z = \sum_{i \in N} c_i y_i \ + \ VOT \ \times \sum_{i \in N} w(i) \tag{22}$$

### 4.3 Discussions on station-based waiting time estimation and feedback loop

**Waiting Time Estimation:** Given a peak hour demand and capacity at EV charging stations, and the resulting demand-to-capacity (D/C) ratio. Using a queueing-based volume-delay function (QVDF) Zhou et al., 2022), the average waiting time $W_{avg}$ for EV charging can be expressed as a function of the D/C ratio: $W_{avg} = f\left(\frac{D}{C}\right)$, where $f$ is a derived relationship based on fluid or stochastic queuing models, reflecting the system's waiting time characteristics.

Building on this, it is possible to aggregate the station-by-station delay, $W_{avg}$, alongside the number of EV charging stations for a high-level regional evaluation. One could devise a macroscopic delay or cost function contingent on the number of EV charging stations at the regional level, considering spatial variations. This would give rise to an analytical representation of the inconvenience cost as a function of demand $A(t)$ versus the cumulative supply $D(t)$. Such a comprehensive analytical approach merits in-depth exploration in future research.

One of the most challenging questions in our research is how to fully consider congestion effects associated with soft capacity constraints. We can see adding a new set of multipliers (**Eq. 22**) can better capture various real-world constraints but could further complicate the proposed solution framework. We will examine the related algorithmic and computational performance issues for in our future study.

**Feedback Metrics:** To be precise in terms of what feedback metrics we employ, it is crucial to delve deeper into the key determinants of user preferences and behaviors related to the EV charging network. First and foremost, the proximity and accessibility of the charging infrastructure play pivotal roles. Here, proximity does not just refer to a straightforward linear distance but considers the detour distance that a user must undertake to access a charging station. On the other hand, accessibility is a broader metric, encompassing not just the ease of reaching a charging point, but also the expected waiting time during peak hours.

Beyond these spatial and temporal metrics, the overarching inconvenience cost combines several factors, encompassing both direct and indirect expenses users might incur when using the EV network. As a first step in refining our feedback, numerical quantifications of proximity and accessibility can be effectively mapped using space-time network models and leveraging the principles of time geography frameworks. This would provide a more holistic understanding of user experiences and challenges, paving the way for strategic enhancements in network planning.

## 5. NUMERICAL RESULTS
We calibrate the market demand and supply of EV charging stations using real-world EVs and charging station data, incorporating models M1, M2, and M3. With the calibrated supply of EV charging stations as a foundation, we employ the medium-term grid-based charging station planning approach (M4) within a multi-level network framework.

### 5.1 Calibration of EV Charging Supply Curve

This study utilizes public EVSE charging station data sourced from the Alternative Fuel Data Center (AFDC), which is captured by the National Renewable Energy Laboratory (NREL). The data encompass Level 1-3 charging station locations, gathered through service providers and community reporting (Brown et al., 2021). For this research, station data was retrieved for the Chicago area, revealing a total of 1328 charging stations, equipped with 3315 charging outlets. Algorithm 2 depicts the process and algorithm for calibrating EV charging station supply curve parameters.

---

**Algorithm 2: Process and Algorithm for Calculating EV Charging Supply Curve Parameters**

Step 1: Utilizing available demand data (e.g., existing EVs), calculate parameters p, q, and m using the Nonlinear Least Squares Method (NLSM). The purpose of the NLSM algorithm is to minimize $\sum_{i=1}^{N}(\hat{f}(x^{(i)}, \boldsymbol{\theta}) - y^{(i)})^2$, function form $\hat{f}(x^{(i)}, \boldsymbol{\theta})$ is parameterized by parameters $\theta_1, \ldots \theta_{parameters}, and\ (x^{(i)}, y^{(i)}), i \in N$ are data points. The minimization is over the parameters $\boldsymbol{\theta}$.

The equation describing the innovation demand curve is shown below:

$$f(t) = \lambda^I(t) = \frac{d(F(t))}{dt} = m\frac{p(p+q)^2 e^{t(p+q)}}{(pe^{t(p+q)}+q)^2}$$

$$F(t) = m\frac{1-e^{-(p+q)t}}{1+\frac{q}{p}e^{-(p+q)t}}, m \text{ is the total number of innovation users (total number of EV on market)}$$

Step 2: With the parameters p and q established, compute $t_{peak}$, the time corresponding to the peak value of $\lambda^I(t)$, and the inflection point $t_{ip}$ of $\lambda^I(t)$. Utilizing the parameters p, q, and $t_{peak}$, determine $t_0$ from results of NLSM. The demand curve and supply curve shared the same starting time of $t_0$.

$$t_{peak} = \frac{1}{p+q}\ln(\frac{q}{p}), \ t_{ip} = t_{peak} \pm \frac{\ln(2+\sqrt{3})}{p+q}$$

Step 3: According to Assumption 1, the net inconvenience $\pi^I(t)$ follows the polynomial distribution, by given or calculated $\rho^I$ and $t_2$, derive function forms $\pi^I(t)$ and $Q_t$. ($\rho^I$ and $t_2$ can receive predefined values from decision-makers, or users map opt to leverage existing supply data, such as the number of EV charging stations, to calculate $\rho^I$ and $t_2$ through the utilization of NLSM.)

Step 4: The function form of the supply curve can be derived using the equation: $\mu^I(t) = \lambda^I(t) - \pi^I(t) = \frac{p(p+q)^2 e^{t(p+q)}}{(pe^{t(p+q)}+q)^2} - \rho^I(t-t_0)(t-t_2)$, the cumulative supply curve $D^I(t)$ can also be calculated.

Considering the supply curve derived in **Eq. (16)**, it becomes evident that the current number of EVs within the research area (as depicted in Fig. 10a) serves as a reliable indicator of the actual demand for electric vehicles. To effectively calibrate the parameters in the supply curve, $\mu^I(t)$ hinges upon three significant elements: firstly, the coefficient of innovation denoted by $p$ plays a vital role in capturing the rate at which innovations and technologies are adopted by consumers, thereby influencing the overall growth of EV adoption. Secondly, the coefficient of imitation, represented by $q$ is instrumental in modeling the impact of social influence and the tendency of individuals to adopt new technologies based on the behavior of others. Lastly, the coefficient $\rho^I$ represents the net flow between the arrival and dispatch rates, reflecting the dynamic interactions between EV charging station arrival and departure rates over time. (Lekvall et al., 1973; Sultan et al., 1990; Mahajan et al., 1995) offer valuable insights into the range of parameters for $p$ and $q$ which represent the adoption of emerging technologies in the market. By calibrating these parameters, the supply curve can be tailored to accurately represent the interplay between demand and supply in the context of EV adoption, aiding in the formulation of informed strategies for the further proliferation of electric vehicles in the research area.

By combining Step 1 in Algorithm 2 and existing demand data in Fig. 10(a), the $p$, $q$, $t_{peak}$ and $t_0$ can be calculated as follows: $p = 0.005$, $q = 0.5$, $t_{peak} = 2026$ (year), $t_0 = 2012$ (year). The coefficient $\rho^I$ directly influences the net inconvenience upon calibrating $\mu^I(t)$ using real-world charging stations data collected for the Chicago area, Fig. 10(b) is the number of charging stations and charging ports in the Chicago area. Having p, q and $t_0$, as explained in Step 4 in Algorithm 2, $\rho^I$ and $t_2$ are calibrated as follows: $\rho^I = -2.34297 * 10^{-4}$, $t_2 = 2030$, $t_3 = 2045$. For an in-depth exploration of the precise derivation and calculations, please consult Luo et al. (2023).

The presented numeric values represent a snapshot of the situation based on a synthesis of historical data and predictive algorithms. However, the model inherently acknowledges the variability of real-world conditions and the fact that demand and supply dynamics can undergo shifts as the years unfold. By incorporating calibration mechanisms such as like $p$, $q$, $\rho^I$, $t_2$, and $t_3$, the model can continually adjust its predictions to match the evolving landscape of charging station demand, ensuring that its insights remain dependable and valuable for decision-makers in the face of changing circumstances.

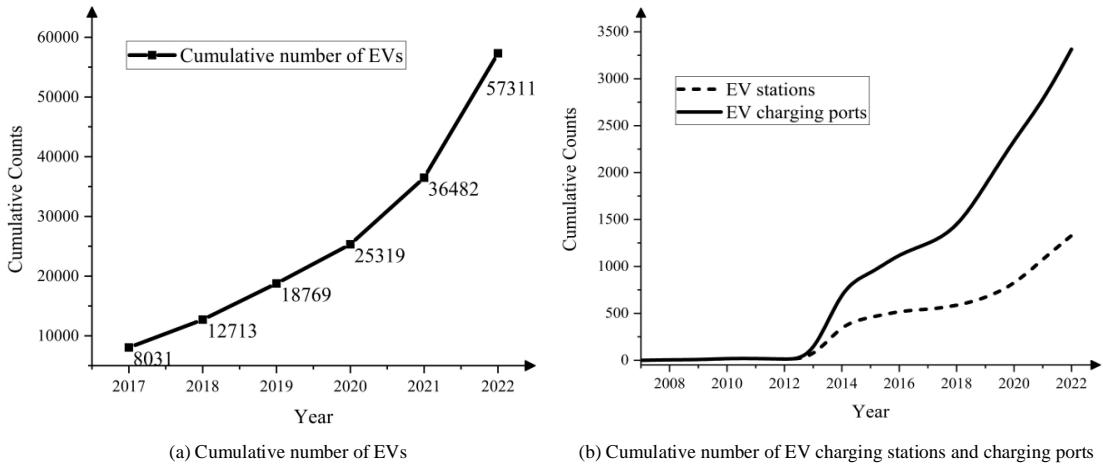

(a) Cumulative number of EVs

(b) Cumulative number of EV charging stations and charging ports

**Fig**. 10 Cumulative number of EVs, Charging Stations and Charging Ports in the Chicago Area

**Fig. 11** showcases the final predictions for the supply of EV charging stations. From the market demand for EV charging stations, the peak demand year is 2026 (Fig. 11a) and the demand will come to its maximum size in the year 2038 (Fig. 11b). EV changing station supply will follow the same peak demand rate in 2026 (Fig. 11e), while the supply will meet its maximum market size at the year of 2045 (Fig. 11f), which is 7 years behind the demand of the market.

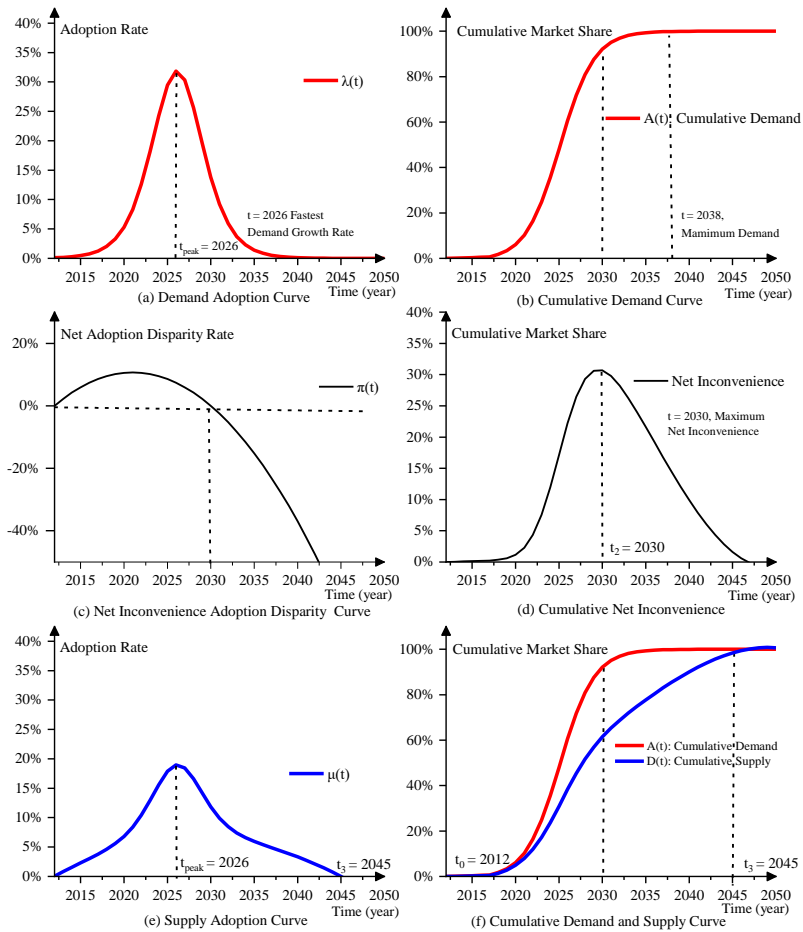

**Fig. 11** Computed Demand and Supply of Electric Vehicle Charging Stations in the Chicago Region

The difference between the cumulative EVs in Fig. 10 and percentage-based representation in Fig. 11 is apparent, arising from differing data presentation approaches. While Fig. 10 displays the actual count of EV charging stations over years, Fig. 11 employs a percentage-based portrayal of cumulative supply for a comprehensive view of market dynamics. The omission of precise station counts, despite calculated supply, is to emphasize evolving dynamics influenced by numerous factors. To maintain precision, we use percentage-based illustrations, capturing relative demand-supply imbalances regardless of numerical shifts. This methodology provides holistic insights into changing market dynamics.

## 5.2 Grid-based EV Charging Location Selection: Network Planning Across Different Resolutions

The study primarily focuses on an aggregated representation of the Chicago sketch network. To generate driving paths within this network, the dynamic traffic assignment simulation tool DTALite (Zhou and Taylor, 2014; Tong et al., 2019; Xiong et al., 2018, 2020) employs the column generation algorithm, effectively generating system-wide driving paths. From these driving paths, specific EV driving paths are selected based on different EV driving distances (MirHassani et al., 2013; Zhou et al 2022; Lu et al., 2022, 2023). We refine the selection process by focusing on two specific round-trips EV driving ranges: 150 km and 250 km in this study. The criteria for selecting these paths are primarily entered around the driving distance. For the 150 km driving range, we choose paths that extend beyond 150 km and disregard those that fall short of this threshold. The rationale behind this approach is that if a path is shorter than 150 km, drivers wouldn't require a charging stop between their origin and destination. Similarly, for the 250 km driving range, we select paths that exceed 250 km while omitting paths that are shorter. By employing this selection process, we concentrate our analysis on paths that align with the driving ranges of interest, ensuring that the charging infrastructure planning is tailored to the specific demands of longer journeys. The selection is predominantly driven by the driving distance criterion, enabling a more targeted investigation into the charging needs of EVs over varying distances.

In the context of EV charging station network planning within a temporal (yearly) horizon, decision-makers need to consider region-based, block-based, and corridor-based deployment strategies. We utilize a three-tier approach, comprising Macro-Meso-Micro levels of EV charging station network planning, to delineate distinct planning stages. At the Macro-Level (10x10 grid), we widen our perspective to encompass entire cities or regions using a generalized grid, facilitating early-stage decisions based on the broader context. As we transition to the Meso-Level (50x50 grid), our focus deepens, allowing us to understand potential EV paths and locations with greater detail, serving as a bridge between macro and micro complexities. Finally, at the Micro-Level (100x100 grid), we delve into meticulous analysis, where each grid cell captures specific city blocks or neighborhoods. This high-resolution evaluation ensures optimal EV station placement, prioritizing accessibility, and operational efficiency. By dividing our approach into Macro, Meso and Micro levels, we ensure thoroughness - from a bird's-eye view down to the intricacies of individual neighborhoods. **Table 3** provides an overview of the selected Origin-Destination (OD) paths within the region, considering driving ranges of both 150 km and 250 km. For an in-depth exploration of the precise calculations and codes, please consult Luo et al. (2023). The location selection algorithm was executed on a CPU model featuring 40 physical cores, 40 logical processors, and 40 threads, utilizing Gurobi Optimizer version 9.1.1 build version 9.1.1rc0 (Linux).

**Table 3** Summary of computational performance in multi-scale networks

| Network Layout (Grid Cells) | Selected # of EV Driving Paths | EV Driving Range R (km) | Capacity Volume | # of Variables: EV Arcs | # of Variables: Candidate Locations | Objective Value (# of stations) | CPU Running Time (s) |
|---|---|---|---|---|---|---|---|
| Macro (10*10) | 19 | 250 | NA | 1,627 | 100 | 2 | 0.10 |
| | 19 | 250 | 5 veh/hour | 1,627 | 100 | 16 | 1 |
| | 1,752 | 150 | NA | 77,999 | 100 | 10 | 8 |
| | 1,752 | 150 | 450 veh/week | 77,999 | 100 | 23 | 32,268 |
| Meso (50*50) | 8 | 250 | NA | 4,693 | 2,500 | 2 | 0.12 |
| | 8 | 250 | 5 veh/hour | 4,693 | 2,500 | 6 | 1 |
| | 8,672 | 150 | NA | 130,113 | 2,500 | 17 | 40 |
| | 8,672 | 150 | 70 veh/week | 130,113 | 2,500 | 45 | 20,378 |

| | 8 | 250 | NA | 4,730 | 10,000 | 3 | 0.14 |
|---|---|---|---|---|---|---|---|
| Micro | 8 | 250 | 5 veh/hour | 4,730 | 10,000 | 5 | 0.36 |
| (100*100) | 607 | 150 | NA | 271,269 | 10,000 | 18 | 358 |
| | 607 | 150 | 80 veh/week | 271,269 | 10,000 | 91 | 124,718 |

As introduced in Eq. (21) and Eq. (22), we employ equal capital cost value over stations. However, we acknowledge that the actual cost of building a station can significantly differ across various locations. While our current paper does not incorporate this future collaboration with decision makers, it is worth noting that our intention is to enhance our approach by gathering region-specific station construction cost data, capacity volume data and comprehensive insights into delay time costs through future collaborations. This prospective collaborative effort aims to align our methodology with real-world parameters, contributing to a more accurate and impactful analysis.

The capacity volume values presented in each instance of **Table 3** are carefully determined through practical considerations and flow requirements. To tailor our approach to different network levels, we employ distinct capacity units: vehicles per hour for micro-level station volume specifications, and vehicles per week for meso-level and macro-level specifications. Drawing from the existing flow volumes on each path, these specific capacity numbers serve as essential lower-bound constraints to maintain the feasibility of the optimization model applied to the network-based station selection process, where two Level-2 charging ports are designed for each selected station. The selection process is grounded in the imperative to strike a balance between accommodating the projected flow of EVs and preserving the network's integrity. If station capacities fall below these predefined boundaries, the optimization model may produce infeasible solutions. Our choice of these capacity values is guided by our aim to pre-empt potential infeasibilities in the optimization process. Lower capacities at stations could require detours and additional paths to accommodate increased demand, although this aspect is not within the scope of our present study but signifies a prospect for exploration in our future research.

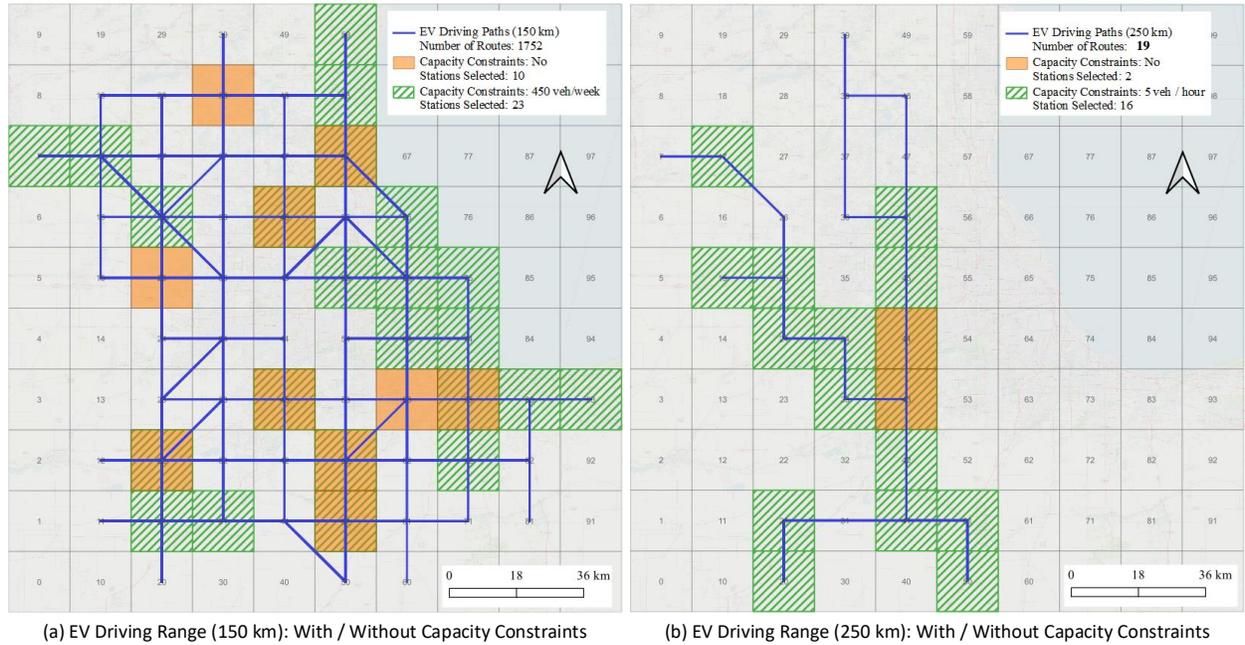

(a) EV Driving Range (150 km): With / Without Capacity Constraints  (b) EV Driving Range (250 km): With / Without Capacity Constraints

**Fig. 12** Macro-Level Selection of EV Charging Locations (10*10 Grid Cells) With / Without Capacity for Candidate Stations

At the macro-level of network planning for EVs with driving ranges of 150 km and 250 km (refer to Fig. 12), the initial analysis encompasses 10 potential locations for the 150 km range and 2 locations for the 250 km range. Initially, our initial model operates under the assumption of unlimited capacity constraints for candidate stations. Subsequently, as we incorporate capacity constraints, we observe significant adjustments that yield more optimal results. Specifically, the findings underscore optimal selections of 23 locations for the 150 km range, accommodating a capacity of 450 vehicles per week, and 16 locations for the 250 km range, with a capacity of 5 vehicles per hour. The identification of these specific sites assumes paramount importance in facilitating judicious decision-making regarding location selection.

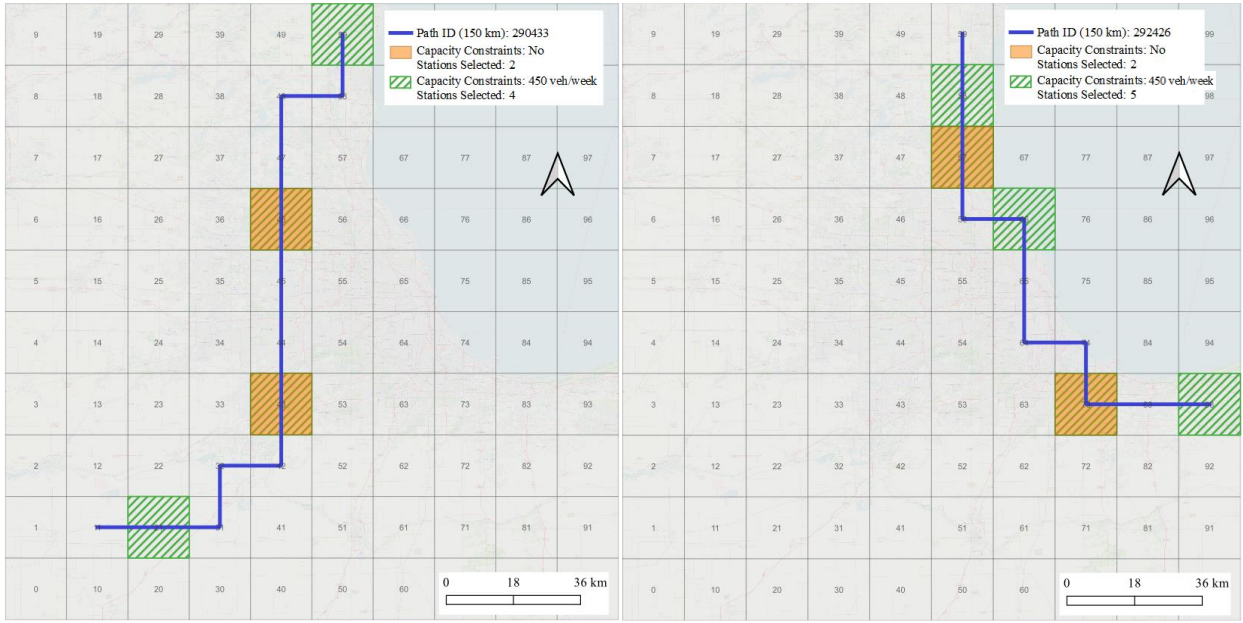

(a) Path ID (150 km): 290433, With/Without Capacity Constraints    (b) Path ID (150 km): 292426, With/Without Capacity Constraints

**Fig. 13** Selection of Single-Driving-Path EV Charging Locations (with and without Capacity Constraints)

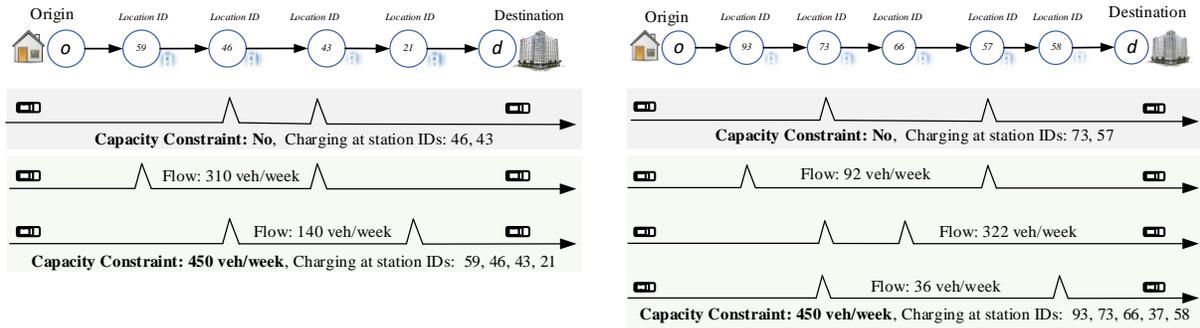

(a) Path ID 290433 (150 km): Segregated Flow with Station-Specific Capacity Constraints    (b) Path ID 292426 (150 km): Segregated Flow with Station-Specific Capacity Constraints

**Fig. 14** Illustration of EV Charging Flows (with and without Capacity Constraints)

Fig. 13 portrays the influence of capacity constraints on the selection of locations for two distinct EV driving routes. Fig. 14 illustrates EV charging flows for two selected paths. For the first route (indexed as 290433) within the 150 km driving range, the absence of capacity constraints results in the selection of locations with ids 46 and 43, necessitating two charging facilities-first at location 46 and then at location 43. However, upon integrating the stipulated capacity of 450 vehicles per week for each station, the optimized choice encompasses locations 59, 43, 46, and 21. This yields a refined allocation of vehicle flows: 22 vehicles recharge at location 46 and subsequently at location 21, whereas 10 vehicles follow a charging sequence of location 59 followed by location 43, culminating in the completion of their journeys on this route.

Similarly, for the second route (indexed as 292426) with a 150 km driving threshold, initial selection in the absence of capacity constraints falls on locations 73 and 57, prompting drivers to charge at location 73 and subsequently at location 57 to finalize their trips. However, when considering the designated capacity of 450 vehicles per week for each prospective station, a total of five stations emerges as optimal choices. Specifically, this involves the assignment of 18 vehicles to the location sequence: start-93-57-destination, 63 flows directed towards the sequence: start-73-66-destination, and 7 flows designated for the sequence: start-73-58-destination. The integration of capacity constraints prompts the selection of a greater number of EV charging locations and fosters a more refined flow assignment strategy, ultimately enriching the overall landscape of network planning.

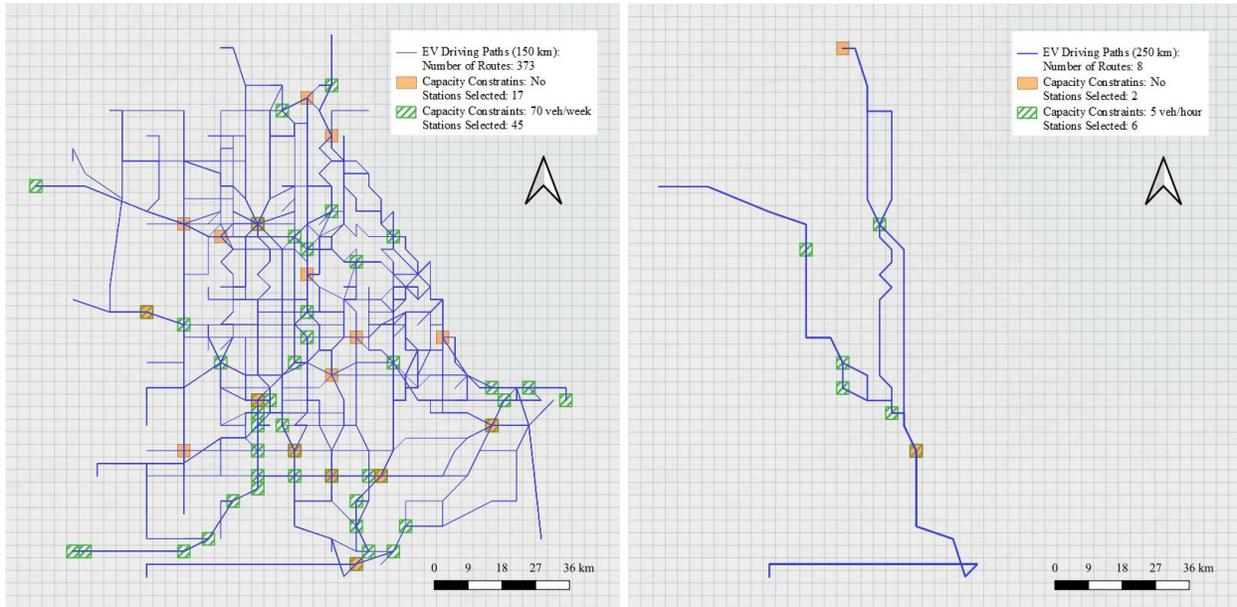

(a) EV Driving Range (150 km): With / Without Capacity Constraints    (b) EV Driving Range (250 km): With / Without Capacity Constraints

**Fig. 15** Meso-Level Selection of EV Charging Locations (50*50 Grid Cells) With/Without Capacity Constraints for Candidate Stations

The meso-level of network planning is characterized by a comprehensive examination of multiple paths and EV driving distances, particularly focusing on distances of 150 km and 250 km. In the exploratory phase, represented by **Fig. 15**, a total of 17 locations are selected for the 150 km driving range, while 2 locations are selected for the 250 km driving range. Notably, at this stage, candidate stations' capacity constraints are not factored into the analysis. To address the crucial aspect of capacity considerations, the optimal locations are identified as 45 for the 150 km driving distance, with a capacity of accommodating 450 vehicles per week, and 6 for the 250 km driving distance, boasting a capacity to serve to 5 vehicles per hour. The strategic selection of these locations holds immense significance in the process of making well-informed and effective decisions when it comes to the placement of EV charging stations.

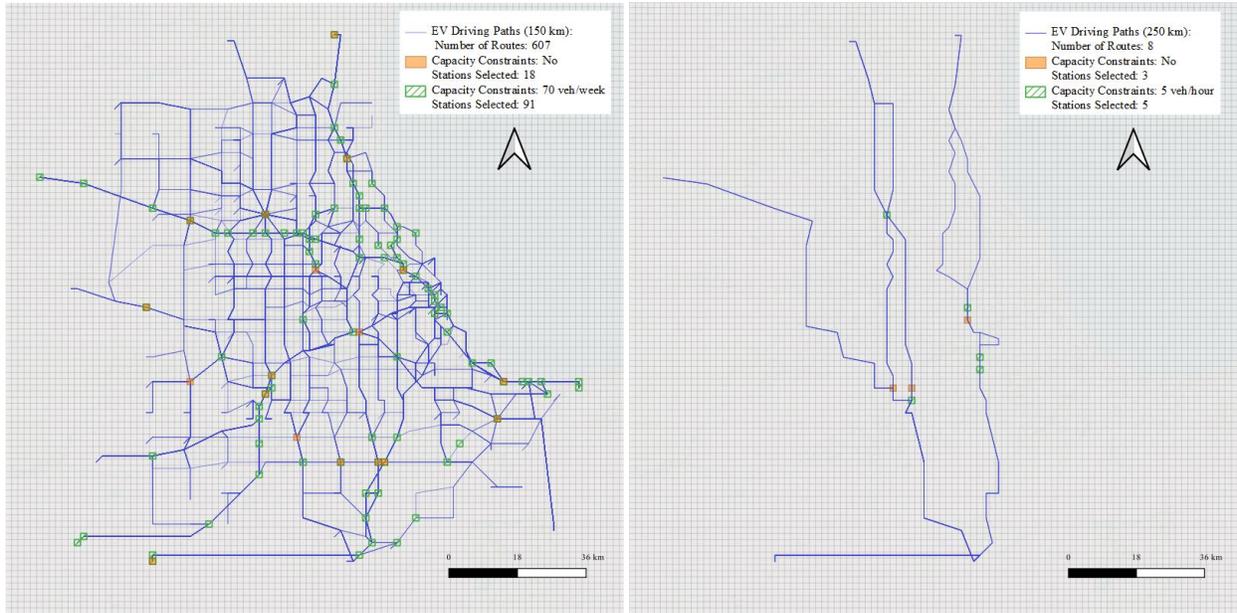

(a) EV Driving Range (150 km): With / Without Capacity Constraints      (b) EV Driving Range (250 km): With / Without Capacity Constraints

**Fig. 16** Micro-Level Selection of EV Charging Locations (100*100 Grid Cells) With / Without Capacity for Candidate Stations

At the micro-level of network planning, the central focus revolves around the meticulous examination of multiple paths and EV driving distances, notably spanning 150 km and 250 km. The initial analysis, as depicted in Fig. 16, encompasses a total of 18 optimal locations for the 150 km driving distance and 3 locations for the 250 km driving distance. Notably, this analysis in Fig. 16 is conducted without taking into consideration the constraints posed by candidate station capacity. To address this pivotal aspect of capacity considerations, this study undertakes a comprehensive exploration of diverse forms of capacities, which includes an evaluation of maximum weekly and daily service volumes. Consequently, the optimal locations are identified as 91 for the 150 km driving distance, boasting a capacity to accommodate 450 vehicles per week, and 5 for the 250 km driving distance, with a capacity to accommodate 5 vehicles per hour. The strategic selection of these locations holds paramount significance in the process of making well-informed and effective decisions regarding the placement of EV charging stations.

## 6. CONCLUSIONS

In conclusion, this paper presents a novel and sophisticated approach for addressing critical aspects of innovation diffusion and supply curve design in electric vehicle charging networks. By employing a continuous-time fluid queue and discrete flow refueling location modeling method, the research offers valuable insights into the dynamic of long-term demand and supply of charging stations, then make a connection to the supply side charging station deployment with optimizing the deployment of electric vehicle charging infrastructure. The integration of the formula to establish an alternate form for the cumulative supply curve enables efficient consideration of multiple deployments among cities within a space-time network model. Moreover, the comprehensive analysis considers both overutilization and underutilization scenarios, making it indispensable for policymakers and stakeholders in developing effective infrastructure plans and conducting informed cost-benefit analyses.

The proposed computational framework empowers decision-makers with accurate estimates and predictions of market adoption rates and potential saturation levels over extended timeframes, facilitating strategic resource planning decisions. Additionally, the extension of the demand-supply model from a single region to multiple regions within a network design framework opens avenues for optimizing EV charging deployment across various cities. These contributions significantly advance the field of transportation research and offer a strong foundation for enhancing the efficiency and sustainability of EV charging networks. As the world continues to embrace electrification and sustainable transportation, this research stands poised to guide the development of intelligent and impactful EV charging solutions for a cleaner and greener future.

In future research, we intend to explore model capacity constraints in more depth, especially their impact on numerical experiments. While this current study provides insights without incorporating numerical experiments related to EV charging capacity constraints, subsequent research will focus on the application of these soft or hard constraints.

Our hypothesis is that models with soft constraints tend to find an equilibrium. They may lean towards adding more stations, which would decrease waiting costs at the expense of higher capital outlay. Conversely, opting for fewer stations might cut capital costs but lead to longer waiting times.

A potential direction is to employ a volume/capacity methodology to estimate the average waiting time in the objective function. This would entail considering queueing behavior with set discharge rates. For instance, a recent study by Zhou et al. (2022) aims to expand upon Newell's fluid queue models to evaluate analytical waiting time functions in scenarios with oversaturated traffic conditions.

In addition to the contributions outlined in this paper, there are several promising avenues for future research in the domain of innovation diffusion and supply curve design within EV charging networks. Firstly, the proposed model could be further validated and refined using extensive real-world data from diverse charging network deployments, allowing for a comprehensive assessment of its predictive accuracy and practical applicability. Secondly, the integration of dynamic pricing strategies could enhance the model's ability to effectively manage network congestion and incentivize optimal charging behaviors. Moreover, the incorporation of user behavior modeling holds the potential to offer deeper insights into adoption patterns, contributing to a more nuanced understanding of the factors influencing EV uptake. Lastly, as the landscape of EV technology rapidly evolves, investigating the impacts of emerging trends such as renewable energy integration, autonomous vehicles, and evolving government policies will be instrumental in maintaining the model's relevance and accuracy in guiding the future development of robust and efficient EV charging networks.

## ACKNOWLEDGEMENTS

The first and last authors are supported by National Science Foundation under grant no. TIP-2303748 titled, "CONNECT: Consortium of Open-source Planning Models for Next-generation Equitable and Efficient Communities and Transportation". We extend our heartfelt appreciation to Fang (Alicia) Tang for her invaluable insights and contributions that greatly enriched this paper. Special thanks to Jason Hawkins for his valuable feedback, which has been instrumental in refining our work.